\documentclass[11pt,a4paper]{article}

\usepackage{setspace}
\usepackage{authblk}
\usepackage{utopia}
\usepackage{enumerate}
\usepackage{amsmath}
\usepackage{graphicx}
\usepackage{amssymb}
\usepackage{amsthm}
\usepackage{stmaryrd}
\usepackage[margin=2.54cm]{geometry}
\usepackage{float}

\makeatletter
\renewcommand*\@fnsymbol[1]{\the#1}
\makeatother

\theoremstyle{plain}
\newtheorem{Theorem}{Theorem}
\newtheorem{Corollary}[Theorem]{Corollary}
\newtheorem{Lemma}[Theorem]{Lemma}

\theoremstyle{definition}
\newtheorem{Definition}[Theorem]{Definition}
\newtheorem{Conjecture}[Theorem]{Conjecture}

\newcommand{\pr}{\text{\textrm Pr}}

\newenvironment{pf}{
    \begin{proof}
        [\textsc{Proof}]
        }
    {\end{proof}
}

\title{{\bfseries Random Information Spread in Networks}}
\author[a,b]{Raymond Lapus \thanks{\tt {jrrlapus24@yahoo.co.nz, rlapus@hs-mittweida.de}}}
\author[b]{Frank Simon \thanks{\tt {simon@hs-mittweida.de}}}
\author[b]{Peter Tittmann \thanks{\tt {peter@hs-mittweida.de}}}

\affil[a]{{\small De La Salle University, Taft Avenue 2401, 1004 Manila, Philippines}}
\affil[b]{{\small Mittweida University of Applied Sciences, Technikumplatz 17, 09648 Mittweida, Germany}}

\begin{document}
    \maketitle

\begin{abstract} 
    Let $G=(V,E)$ be an undirected loopless graph with possible parallel edges
    and $s,t\in V$. Assume that $s$ is labelled at the initial time step and that 
    every labelled vertex copies its labelling to neighbouring vertices 
    along edges with one labelled endpoint independently with probability
    $p$ in one time step. In this paper, we establish the equivalence between the 
    expected $s$-$t$ first arrival time of the above spread process 
    and the notion of the stochastic shortest $s$-$t$ path. Moreover, we
    give a short discussion of analytical results on special graphs including the
    complete graph and $s$-$t$ series-parallel graphs. Finally we propose
    some lower bounds for the expected $s$-$t$ first arrival time.

    \vspace{0.125in}
\end{abstract} 

    \section{Introduction}
    The spreading of information in networks is a random process. 
    Consider an undirected loopless graph that possibly has parallel edges
    $G = (V,E)$. Let $s \in V$ be a chosen vertex that is labelled at time step 0. 
    In addition we assume that that every edge of $G$ that has a labelled endpoint
    independently copies the label to the possibly unlabelled endpoint with the 
    given infection probability $p$ in one time step.

    Some applications are network based models for virus spread
    in epidemiology or the analysis of gossipping across social networks. The
    analysis of malware propagation in computer networks, such as spread of
    computer viruses or worms, represents another important application in
    technical networks. 
    
    The epidemic threshold is an interesting measure that characterises
    the transition from local to global spread. Random graphs and small-world
    networks gained most attention for investigations in this field (see for
    instance \cite{ARBA02:statmechanics} or \cite{NM03:strucfunccplxnetwork}).
    Certain regularity conditions are necessary in order to make the analysis
    of large or infinite graphs possible. An application for the analysis of
    the outbreak of Severe Acute Respiratory Syndrome (SARS) is presented in
    \cite{MLPBNM05:netwrokssars}.

    We deal with the following questions in this paper. First we present the problem of 
    finding the expected $s$-$t$ first arrival time using the generating 
    function approach (see Section \ref{sec:spreadprocess}). Second we establish in Section
    \ref{sec:stochasticShortestPaths} the relationship between the spread process and the
    stochastic shortest $s$-$t$ path problem stated in \cite[Chapter 9]{SD91:netrelalgstrucs}).   
    In Section \ref{sec:spreadingResistance} we introduce the $s$-$t$ spreading resistance of
    a graph. We show that
    this notion is related to Kulkarni's concept \cite{KV86:spathexpdistarclength}
    of the expected length of a shortest $s$-$t$ path, whenever the edge
    lengths are exponentially distributed random variables with intensity $p$.
    Subsequently, we establish two reduction techniques  for calculating the expected first 
    arrival time in $s$-$t$ series-parallel
    graphs (see Section~\ref{sec:reductiontechs}). Due to the computational complexity of the
    spread process, we propose
    in Section \ref{sec:bounds} methods for yielding lower and upper bounds for
    the expected $s$-$t$ first arrival time
    in terms of the $s$-$t$ reliability polynomial and the distance 
    between $s$ and $t$.

    \section{The spreading process\label{section:spreadingprocess}}
    \label{sec:spreadprocess}
    All graphs considered in this paper are finite and loopless, but 
    possibly with parallel edges. In this section, we present how the
    spread process propagates in arbitrary graphs.
    
    \subsection{Assumptions} 
        Let $G=(V,E)$ be a graph and
        $X_k \subseteq V$ be the set of labelled vertices in $G$ at the
        discrete nonnegative integral time step $k$. Assume that $u \in X_k$
        copies the message to an adjacent unlabelled vertex $v$ with
        probability $p_{uv} = 1 - q_{uv}$ along the edges $\{u,v\}$ in $G$. 
        We refer to the numbers  $p_{uv}$ and $q_{uv}$ as the \emph{infection} and 
        \emph{noninfection probability from $u$ to $v$}.

        The spread is performed simultaneously along all edges that link
        labelled vertices to unlabelled vertices within one unit of time. 
        The copying of labels along the edges is assumed to be stochastically independent.
        In addition the message spread is symmetric, that is for all $u,v \in V$,
        $p_{uv} = p_{vu}$. We presuppose that $p_{uv} > 0$ for any $\{u,v\} \in E$,
        whereas $\{u,v\} \notin E$ implies $p_{uv} = 0$.

        Therefore we can describe the spread process by a homogeneous 
        Markov chain $\{ X_k \}_{k \in \mathbb{N}}$, where the state space is a 
        subset of $2^V$. If $G$ is a connected graph then the labelled vertices spread
        until all vertices of $G$ are labelled; that is, $\lim_{k \rightarrow \infty} X_k = V$.

        Let $G=(V,E)$ be a graph and $v \in V$. The \emph{open neighbourhood of $v$},
        denoted by $N(v)$, is the set of vertices that are adjacent to $v$ in $G$. For each
        $A \subseteq V$, the \emph{open neighbourhood of $A$} is defined by
        \begin{align*}
            N(A) := \bigcup_{v \in A} N(v) \setminus A.
        \end{align*}
        Furthermore for $A,B\subseteq V$, the \emph{$A$-$B$ edge set} in $G$, 
        denoted by $(A,B)$, is the set of edges of $G$ with one endpoint in $A$ 
        and the other endpoint in $B$. That is,
		\begin{align*}
			(A,B)&:=\{\{u,v\} \in E \colon u \in A, v \in B\}.
		\end{align*}
        \begin{Definition}
        Let $G=(V,E)$ be a graph and $A,B\subseteq V$ with $A\subseteq B$. The 
        \emph{transition probability} $P_G(A,B)$ is the probability that the vertex 
        subset $B$ becomes infected in one time step if the vertex subset $A$ is infected 
        beforehand. 
        \end{Definition}
        \begin{Theorem}
            Let $G=(V,E)$ be a graph and $A,B \subseteq V$ with $A \subseteq B$.
			   The transition probability $P_G(A,B)$ is then given by
            \begin{align*}
                P_G(A,B) = \prod_{v \in A} \prod_{e \in (\{v\}, N(A) \setminus B)} q_e
                        \prod_{x \in B \setminus A} \left[1 - \prod_{f \in (\{x\},A)} q_f\right].
            \end{align*}
            Moreover if $G$ is a simple graph, then
            \begin{align*}
                P_G(A,B) = \prod_{v \in A} \prod_{w \in N(A) \setminus B} q_{vw}
                    \prod_{x \in B \setminus A} \left[1 - \prod_{y \in N(x) \cap A} q_{xy}\right].
            \end{align*}
        \end{Theorem}

\begin{figure}
        		\begin{center}
        		    \includegraphics{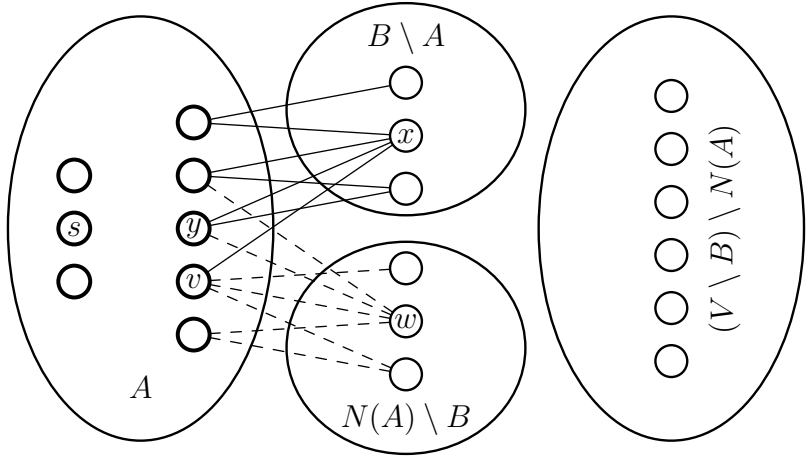}
        		\end{center}
            \caption{Transition probability of the spread process from $A$ to $B$ in $G$
            }
            \label{fig:transition}
        \end{figure}

        \begin{pf}
            Figure~\ref{fig:transition} explains the derivation of $P_G(A,B)$ in a
            simple graph $G$. There must be no label propagation along edges from
            $A$ to $N(A) \setminus B$ that are represented as broken lines. This
            gives the first double product. The remaining products correspond to
            transmissions along solid drawn edges. For each vertex ~$x \in B \setminus A$,
            there must be at least one edges that transports the label which gives
            the term in brackets.
        \end{pf}

        \noindent At this point, we denote the transition probability from $A$ to $B$ 
        in the graph $G$ as $P(A,B)$, whenever $G$ is known from the context.
        
        \begin{Corollary}
            \label{cor:transprob}
            Let $G=(V,E)$ be a graph and $A, B \subseteq V$ with $A \subseteq B$.
            Furthermore, assume that all infection probabilities are identical; 
            that is, $p_{uv} = p = 1-q$ for all $\{u,v\} \in E$. The 
            transition probability from $A$ to $B$ in $G$ is then given by
            
            $~$ \hfill %
            $\displaystyle P(A,B) = q^{|(A,N(A) \setminus B)|}
                 \prod_{x \in B \setminus A} \left[ 1 - q^{|(\{x\},A)|} \right]$.
            \hfill $\blacksquare$
        \end{Corollary}
    
    \noindent Note that $P(A,B)$ can be only nonzero if $B \subseteq N(A) \cup A$ holds.

    \subsection{First arrival times}
        We first describe some general properties of the spread process. In
        order to simplify the presentation, we now assume that all infection
        probabilities are identical.

        \begin{Definition}
            Let $G = (V,E)$ be a graph and $X \subseteq V$ be a nonempty vertex subset. 
            The graph $G_X$ is obtained by \emph{merging} the vertex subset $X$ 
            in $G$. The \emph{merged vertex} in $G_X$ is denoted by $X'$, where $G_X$ contains
            an edge $\{X',y\}$ for each edge $\{x,y\}\in E$ in $G$ with $x\in X$ and 
            $y\in V\setminus X$.
        \end{Definition}
        
        \begin{Lemma}
            Let $G = (V,E)$ be a graph and $A,X \subseteq V$. 
            If $A'= (A \setminus X) \cup \{X'\}$ is set, then \linebreak 
            $P_G(X,A) = P_{G_X}(\{X'\},A')$ holds.
        \end{Lemma}

        \begin{pf}
            Let $N_G(X)$ and $N_{G_X}(\{X'\})$ be the corresponding open neighbourhood 
            of $X$ in $G$ and $\{X'\}$ in $G_X$,  
            then $|(X,N(X) \setminus A)| = |(\{X'\}, N(\{X'\}) \setminus A)|$ holds.
            The merging operation provides $|(X,\{a\})|_G = |(\{X'\},\{a\})|_{G_X}$ for all
            $a \in A \setminus X$. Consequently, the exponents of $q$ from Corollary
            \ref{cor:transprob} remain unchanged while transforming $G$ into $G_X$.
        \end{pf}

        Thus it is advantageous to work with $G_{X}$ instead of $G$ since the
        number of vertices is reduced, but the graph $G_{X}$ may not be simple even in the
        case that $G$ is a simple graph. However, $G_{X}$ can be transformed into a simple
        graph by replacing all parallel edges joining a pair of vertices $u,v \in V$
        by just one edge $\{u,v\}$. That is, two parallel edges namely, $e = \{u,v\}$ and
        $f = \{u,v\}$ with infection probabilities $p_e$ and $p_f$ are replaced by a
        new edge $g$ with infection probability $p_g = p_e + p_f - p_e p_f$. In this
        case we obtain the more general model including edges with different probabilities.

        \begin{Definition}
        Let $G=(V,E)$ be a graph and $A,B\subseteq V$. Denote by 
        $Z_{AB}$ the random variable of the first time that all vertices in the set $B$ are 
        infected if the vertex set $A$ is infected at the beginning of the spread 
        process.
        \end{Definition}

        The following lemma describes the transition from one state to
        another in terms of a Markov chain. This result will be utilised
        to establish Theorem \ref{thm:fundamentalRecursion}. In the following 
        $\subsetneq$ denotes the proper subset inclusion.

        \begin{Lemma}
            \label{lem:transitionMC}
            Let $G=(V,E)$ be a graph and $A,B \subseteq V$ such that $A \cap B \subsetneq B$, 
            then
            \begin{align}
                \mathrm{Pr}(\{Z_{AB} = n\}) =
                    \sum_{C \supseteq A} P(A,C) ~\mathrm{Pr}(\{Z_{CB} = n-1\})
\label{eqn:transitionMC}
            \end{align}
            holds for all $n \geq 1$.
        \end{Lemma}

        \begin{pf}
            The condition $A \cap B \subsetneq B$ implies that there is at least one vertex
            $t \in B$ that is not yet infected, when the spread process reaches the vertex
            subset $A$. The probability that all vertices in $B$ are infected for the first
            time after $n$ time steps can be divided into disjoint events by considering one
            further time step. In one time step a superset $C$ of $A$ can be infected, where
            the transition probabilities are given by the $P(A,C)$. In order to infect $B$
            for the first time after $n$ time steps, it is necessary to reach $B$ from $C$ in
            exactly $n-1$ time steps for the first time to account for the already elapsed
            extra time step.
        \end{pf}


        \begin{Definition}
            Let $G=(V,E)$ be a graph and $A,B \subseteq V$.  
            The \emph{ordinary generating function 
            of the $A$-$B$ first arrival times} is defined as
            \begin{align*}
                \Phi_{AB}(G,z) = \sum_{n \geq 0} \mathrm{Pr}(\{Z_{AB} = n\}) z^n.
            \end{align*}
        \end{Definition}

        \noindent It clearly follows from the definition of $Z_{AB}$ that 
        $\Phi_{AB}(G,z) = 1$ holds for the case of $A \supseteq B$.

        \begin{Theorem}
            \label{thm:fundamentalRecursion}
            The ordinary generating function
            of the $A$-$B$ first arrival times satisfies the recurrence relation
            \begin{align*}
                \Phi_{AB}(G,z) = \sum_{C\supseteq A} z ~P(A,C) ~\Phi_{CB}(G,z)
            \end{align*}
            for all $A,B \subseteq V$ with $A \cap B \subsetneq B$, where the initial
            conditions $\Phi_{CB}(G,z) = 1$ hold for all $C \supseteq B$.
        \end{Theorem}

        \begin{pf}
            The statement follows directly from Lemma~\ref{lem:transitionMC}. 
            Note that the assumption of $A \cap B \subsetneq B$ implies
            $\mathrm{Pr}(\{Z_{AB} = 0\}) = 0$, as $B$ cannot be reached from 
            $A$ in zero time steps.
        \end{pf}

        Observe that the nice triangular structure of the system of linear equations in
        Theorem \ref{thm:fundamentalRecursion}, caused by $P(A,B)=0$ for all pairs
        $(A,B)$, whenever $A$ is not a subset of $B$, permits the practical solution of systems
        with more than $10^{5}$ equations. It is now possible to deduce a recurrence
        relation for the expectation of the random variable $Z_{AB}$ by using the 
        generating function for the $A$-$B$ first arrival time probabilities. A similiar 
        approach for the later discussed exponential model is presented in 
        \cite{KV86:spathexpdistarclength}.

        \begin{Corollary}
            \label{cor:momentreceqn}
            Let $G = (V,E)$ be a graph and $A,B \subseteq V$ with
            $A \cap B \subsetneq B$. The expectation of the random variable
            $Z_{AB}$, denoted by $T_{AB}(G)$, obeys the recurrence equation 
            \begin{align}
            		\label{eqn:momentreceqn}
                T_{AB}(G) = \mathbb{E}[Z_{AB}] =
                		\sum_{C \supseteq A} P(A,C) ~\left[T_{CB}(G) + 1\right]
            \end{align}
            with the initial conditions $T_{AB}(G) = 1$, for any $A \supseteq B$.
             \hfill $\blacksquare$
        \end{Corollary}
        In the case that the sets $A$ and $B$ are singletons say $A=\{s\}$, $B=\{t\}$,
        we will write $T_{st}(G)$ instead of $T_{\{s\}\{t\}}(G)$. This notation is also 
        used when
        only one of the two sets is a singleton, e.g. $T_{At}(G)$ means $T_{A\{t\}}(G)$. The 
        same notation also applies to the random variable $Z_{AB}(G)$ and the later defined spreading
        resistance $\rho_{At}(G)$.

    \section{Stochastic shortest paths}
    \label{sec:stochasticShortestPaths}
    In this section, we recast the problem of the spread process
    in a stochastic shortest path problem. A nice introduction to solve
    stochastic shortest path problems in an algebraic way in
    the case of discrete arc lengths assuming finitely many values
    is found in \cite[Chapter 9]{SD91:netrelalgstrucs}.

    \begin{Definition}
    Let $G=(V,E)$ be a graph and $s,t\in V$. Furthermore let $D_{st}$ be the random variable of
    the length of a shortest $s$-$t$ path in $G$, if the lengths $L(e)$ of the edges $e\in E$ 
    in $G$ are independent geometric random variables with parameter $p$, 
    where $p$ equals the success (infection) probability of the spread process.
    \end{Definition} 

    \begin{Theorem}
        Let $G = (V,E)$ be a graph and $s,t \in V$, then $D_{st}=Z_{st}$.
    \end{Theorem}

    \begin{pf}
    		For all $e \in E$, we let $L(e)$ be 
    		independent random variables representing the edge lengths, that are 
    		drawn from a geometric distribution with success probability
    		$p = 1-q$. That is for all $n \geq 1$,
    		\begin{align*}
    				\mathrm{Pr}(\{L(e) = n\}) = (1-q) q^{n-1}.
    		\end{align*}
			Now for all $A \subseteq V$ the random variable $D_{At}$ denotes 
			the shortest length among all possible lengths of shortest $u$-$t$ paths in $G$
			with $u\in A$, that is,
    				$D_{At} = \mathrm{min} \{D_{ut} \colon u \in A\}$.
    		Note that whenever the target vertex $t$ satisfies $t \in A$, we have 
    		$\mathrm{Pr}(\{D_{At} = 0\})=1$.
    
				Now for every $A\subseteq V$ and $v \in N(A)$, let the random variable $L(A,v)$ be 
				the smallest edge length among all lengths of edges contained in 
				$(A,\{v\})$. That is,
				\begin{align*}
						L(A,v) = \mathrm{min} \{L(e) \colon e \in (A, \{v\})\}.
				\end{align*}
				We observe that
				\begin{align*}
						\mathrm{Pr}(\{L(A,v) = 1\}) = 1 - q^{|(A,\{v\})|}.
				\end{align*}
				With the notions of the length $L(e)$ of an edge $e$ and $L(A,v)$, 
				we can restate the random variable $D_{At}$ as
				\begin{align*}
						D_{At} &= \min \{D_{vt} + L(\{u,v\}) \colon \{u,v\} \in (A,N(A))\} \\
					  &= \min \{D_{vt} + L(A,v) \colon v \in N(A)\},
				\end{align*}
				whenever $A\subseteq V\setminus\{t\}$ holds.
				By defining the event 
				\begin{align*}
						\Lambda(A,B) &:= \{L(A,v) = 1 \colon v \in B \} \cap 
								\{L(A,v) > 1 \colon v \in N(A) \setminus B\}
				\end{align*}
				for each $A\subseteq V$ and $B \subseteq N(A)$, we obtain
				\begin{align*}
						\mathrm{Pr}(\Lambda(A,B)) =
								\prod_{v \in B} \mathrm{Pr}(\{L(A,v) = 1\})
										\prod_{v \in N(A) \setminus B} \mathrm{Pr}(\{L(A,v) > 1\}) 
								= P(A,A\cup B),
				\end{align*}
				which coincides to the transition probability of Section~\ref{section:spreadingprocess}.
				The application of the law of total probability on the event $\{D_{At}=n\}$ yields
				\begin{align*}
						\mathrm{Pr}(\{D_{At}=n\}) &=
								\sum_{B \subseteq N(A)} P(A,A \cup B) 
										~\mathrm{Pr}(\{D_{At}=n\} ~\vert~ \Lambda(A,B)).
				\end{align*}
				For each $A\subseteq V\setminus\{t\}$ and $B\subseteq N(A)$, we have
				\begin{align*}
						\pr(\{D_{At}=n\} ~\vert ~\Lambda(A,B)) &=
								\pr \left(\left\{\min_{v\in N(A)} D_{vt} + L(A,v) = n\right\} 
										\middle\vert ~\Lambda(A,B)\right).
				\end{align*}
				Now the conditional probability on the right hand side can be restated by 
				reducing for all $v\in N(A)$ the 
				realisation of the random variables $L(A,v)$ by one and asking for the probability 
				of a shortest $A$-$t$ path of length $n-1$.  That is,
				\begin{align*}
					\pr(\{D_{At}=n\}\vert ~\Lambda(A,B))&=
						\pr\left(\left\{\min_{v\in N(A)} D_{vt} + L(A,v) = n-1\right\} 
										\middle\vert~\Lambda'(A,B)\right)
				\end{align*}
				with 
				\begin{align*}
					\Lambda'(A,B)=\{L(A,v)=0:v\in B\}\cap \{L(A,v)\ge 1:v\in N(A)\setminus B\}.
				\end{align*}
				Now all $v\in B$ with $L(A,v)=0$ can be included in the set $A$, as edge lengths of 
				size zero do not contribute to the length of a shortest $A$-$t$ path,
				whereas for each $v\in N(A)\setminus B$ the condition $L(A,v)\ge 1$ is always 
				satisfied. Thus the condition $\Lambda'(A,B)$ can be omitted and
				\begin{align*}
						\pr(\{D_{At}=n\} ~\vert ~\Lambda(A,B))
						&= \mathrm{Pr}(\{D_{A \cup B, t} = n-1\})
				\end{align*} 
				follows. Hence, the recurrence relation
				\begin{align*}
						\mathrm{Pr}(\{D_{At} = n\}) &= \sum_{B\subseteq N(A)} 
								P(A,A \cup B) ~\mathrm{Pr}(\{D_{A \cup B, t} = n-1\}) \\
						&= \sum_{C \supseteq A} P(A,C) ~\mathrm{Pr}(\{D_{Ct} = n-1\}),
				\end{align*}
				is valid for all $A \subseteq V \setminus \{t\}$ with initial conditions 
				$\pr(D_{Ct}=0)=1$ for all $C\subseteq V$ with $t\in C$. This recurrence relation 
				coincides with \eqref{eqn:transitionMC}, whenever $Z_{At}$ is replaced 
				by $D_{At}$. Hence, $D_{At} = Z_{At}$ is valid for all $A \subseteq V$ and so
				$D_{st} = Z_{st}$ follows.
    \end{pf}    

    \noindent An important consequence of the above correspondence is the symmetry
    of the $s$-$t$ first arrival time probabilities, which is only obvious in the view of the
    stochastic shortest path formulation of the spread process.

    \begin{Corollary}
        \label{cor:symmetry}
        Let $G=(V,E)$ be a graph and $s,t\in V$ then $Z_{st}=Z_{ts}$.
        \hfill $\blacksquare$
    \end{Corollary}

    \section{Spreading resistance}
    \label{sec:spreadingResistance}
    The spread process is connected to a continuous time Markov chain, which was
    already examined by Kulkarni and Corea in \cite{KV86:spathexpdistarclength,
    CGKV93:spathdiscdistarclength}.

    \begin{Definition}
        Let $G=(V,E)$ be a connected graph with $A \subseteq V$ and $t\in V$.
        Let $T_{At}(G)$ be the expected $A$-$t$ first arrival time in $G$. The
        number $\rho_{At}(G)$ defined by
        \begin{align*}
            \rho_{At}(G) = \lim_{q \to 1} (1-q) T_{At}(G),
        \end{align*}
        whenever this limit exists, is called the \emph{$A$-$t$ spreading resistance of $G$}.
    \end{Definition}

    \begin{Theorem}
        \label{thm:spreadresrecform}
        Let $G=(V,E)$ be a graph then the $A$-$t$ spreading resistance obeys the 
        recurrence relation
        \begin{align*}
            \rho_{At}(G) = \frac{1}{|(A,N(A))|}
                \left[ 1 + \sum_{c \in N(A)} |(A,\{c\})| ~\rho_{A \cup \{c\},t}(G) \right]
        \end{align*}
        for all $A \subseteq V$ and $t\in V$ with $(A,N(A))\ne\emptyset$ and initial conditions 
        $\rho_{Ct}(G) = 0$ for all $C\subseteq V$ satisfying $t\in C$.
    \end{Theorem}

    \begin{pf}
        Observe that the system of linear equations for the expected $A$-$t$ 
        first arrival times that is obtained from \eqref{eqn:momentreceqn} with 
        the initial condition $T_{Ct}(G) = 0$, for all $C \subseteq V$ and 
        $t \in C$ can be restated as
        \begin{align*}
            T_{At}(G) = \frac{1}{1-P(A,A)}
                \left[ 1 + \sum_{C \supsetneq A} P(A,C)~T_{Ct}(G) \right].
        \end{align*}
        As a consequence, the limit value of $(1-q)~T_{At}(G)$ yields to
        \begin{align}
            \label{eqn:limitEqn}
            \lim_{q \to 1}(1-q)~T_{At}(G) =
                \lim_{q\to 1} \frac{1-q}{1-P(A,A)}
                \left[ 1 + \sum_{C \supsetneq A} P(A,C)~T_{Ct}(G) \right]
        \end{align}
        as $q$ tends to $1$. Now consider the two limits from the right-hand side
        of \eqref{eqn:limitEqn}
        \begin{align*}
            L_1 = \lim_{q \to 1}\frac{1-q}{1-P(A,A)} \quad \text{and} \quad
            L_2 = \lim_{q \to 1}\frac{P(A,C)}{1-P(A,A)}.
        \end{align*}
        By applying Corollary \ref{cor:transprob} and L'H\^{o}pital's rule to the
        right-hand side of $L_1$, one finds
        \begin{align*}
            L_1 = \lim_{q \to 1}\frac{1-q}{1-q^{|(A,N(A))|}} = \frac{1}{|(A,N(A))|}
        \end{align*}
        as $|(A, N(A))| \geq 1$. Similarly,
        \begin{align}
            \nonumber
            L_2 & = \lim_{q \to 1} \frac{q^{|(A,N(A)\setminus C)|}}{1-q^{|(A, N(A))|}}
                \prod_{c\in C\setminus A}\left[1-q^{|(A,\{c\})|}\right] \\
            & 
            \label{eqn:limit2simplified}
            =  \sum_{c ~\in ~C\setminus A} \frac{|(A, \{c\})|}{|(A,N(A))|}
                \lim_{q\rightarrow 1} 
                \prod_{\substack{v ~\in ~C\setminus A\\ v \neq c}}
                    \left[1-q^{|(A,\{v\})|}\right].
        \end{align}
        Now suppose $|C\setminus A| \geq 2$. Then for every $c ~\in ~C \setminus A$
        there is a $v \neq c$ that contributes the factor $1-q^{|(A,\{v\})|}$ to the
        product in \eqref{eqn:limit2simplified}. Hence the limit vanishes for every
        $c ~\in ~C \setminus A$. In the case of $C = A \cup \{c\}$ with $c\in N(A)$ one readily finds
        \begin{align*}
            L_2 &= \frac{|(A, \{c\})|}{|(A,N(A))|}.
        \end{align*}
        Now assume that $\rho_{A \cup \{c\},t}(G)$ exists for each $c \in N(A)$ 
        then \eqref{eqn:limitEqn} becomes
        \begin{align*}
            \rho_{At}(G) &= \frac{1}{|(A,N(A))|}\left[1 + \sum_{c \in N(A)}
                |(A, \{c\})| ~\rho_{A \cup \{c\},t}(G) \right]. \hfill \qedhere
        \end{align*}
    \end{pf}

    \begin{Corollary}
        Let $G=(V,E)$ be a connected graph with $A\subseteq V$ and $t\in V$, then the limit
        \begin{align*}
            \rho_{At}(G)&=\lim_{q\rightarrow 1} (1-q)T_{At}(G)
        \end{align*}
        always exists. \hfill $\blacksquare$
    \end{Corollary}

    \begin{Definition}
        Let $G = (V,E)$ be a graph. The \emph{exponential spreading model}
        in $G$ is a spread process such that every edge $e \in E$ is weighted
        with $\omega(e)$ and $\{\omega(e) \colon e \in E\}$ is a collection of
        independent exponential random variables with intensity $p$.
    \end{Definition}

    The following result, due to Kulkarni, shows the connection between
    $\rho_{At}(G)$ and the exponential spreading model in $G$.

    \begin{Theorem}[Kulkarni \cite{KV86:spathexpdistarclength}]
        \label{thm:exponentialModel}
        Let $G=(V,E)$ be a graph with $A \subseteq V$, $t \in V$, and $(A,N(A))\ne\emptyset$. If
        $\tau_{At}(G)$ is the expected length of a shortest $A$-$t$ path
        of the exponential spreading model in $G$, then $\tau_{At}(G)$
        satisfies the recurrence relation
        \begin{align*}
            \tau_{At}(G) = \frac{1}{|(A, N(A))|}
                \left[ \frac{1}{p} + \sum_{c \in N(A)} |(A, \{c\})|
                ~\tau_{A \cup \{c\},t}(G) \right]
        \end{align*}
        with the initial condition $\tau_{Ct}(G) = 0$ for all $C\subseteq V$ with $t\in C$. \hfill $\blacksquare$
    \end{Theorem}

    \noindent The substitution of $(1-q)\tau_{At}(G) = \rho_{At}(G)$ transforms
    Kulkarni's result into the recursive definition of $\rho_{At}(G)$. This
    yields an interpretation of the $A$-$t$ spreading resistance $\rho_{At}(G)$.

    \section{Special graphs}
	\label{sec:specialGraphs}
    This section illustrates the problem of the spread process discussed
    in Sections \ref{sec:spreadprocess} and  \ref{sec:spreadingResistance}
    on graphs with a special structure. The simplest is a path $P_n$ 
    of length $n$. If $s$ and $t$ are the end vertices of $P_n$ and $p$ is
    the infection probability assigned to every edges in $P_n$, then
    $T_{st}(P_n) = n/p$ and $\rho_{st}(P_n) = n$. This result is useful
    when deriving some simple upper bounds as seen in Chapter \ref{sec:bounds}.
    Another trivial example is a tree, which follows from the fact 
    that every two (distinct) vertices in a tree are joined by a unique path 
    of nonzero length.

    \begin{Theorem}
        Let $H(m_1,m_2,\ldots, m_n) = (V,E)$ be the graph consisting
        of $n$ parallel paths of lengths $m_i$ with $1 \leq i \leq n$
        between the vertices $s$ and $t$. Then its $s$-$t$ spreading
        resistance is given by
        \begin{align*}
            \rho_{st}(H(m_1, m_2, \ldots, m_n)) =
                \sum_{j \geq 0} \sum \frac{1}{n^{j+1}}
                    \genfrac(){0cm}{}{j}{i_1, i_2, \ldots, i_n},
        \end{align*}
        where the inner sum runs for all nonnegative integers 
        $i_1, i_2, \ldots, i_n$ such that $j = i_1 + i_2 + \ldots + i_n$ 
        and $i_k < m_k$ holds for all $k$ with $k = 1,2,\ldots,n$.
    \end{Theorem}

    \begin{pf}
        Suppose that for all $j = 1,2,\ldots,m_i$, the weights 
        $\omega(e_{ij})$ of the edges $e_{ij}$ in the $i$-th path of 
        $G = H(m_1,m_2,\ldots,m_n)$ are independent exponential random 
        variables with intensity $p$. Hence the realisation 
        $(\omega(e_{i1}), \omega(e_{i2}), \ldots, \omega(e_{i{m_i}}))$ of 
        waiting times for the edges of the $i$-th path in $G$ corresponds
        to the sum of the random variables $\omega(e_{ij})$, $j = 1,2,\ldots,m_i$.
        This implies that $Y_i = \omega(e_{i1}) + \omega(e_{i2}) + \ldots + 
        \omega(e_{i{m_i}})$ follows an Erlang-$m_i$ distribution with intensity 
        $p$; that is,
        \begin{align*}
            \mathrm{Pr}(\{Y_i \leq t\}) = 1 - \mathrm{e}^{-pt} 
                \sum_{j=0}^{m_i-1} \frac{(pt)^j}{j!},
        \end{align*}
        for each $i = 1,2,\ldots,n$. Now define the random variable
        $Y = \min \{Y_i \colon i = 1,2,\ldots,n\}$, then
        \begin{align*}
            \mathrm{Pr}(\{Y \leq t\}) 
                = 1 - \prod_{i=1}^n \mathrm{Pr}(\{Y_i \geq t\}).
        \end{align*}
        Because $Y$ is a continuous random variable that takes nonnegative
        real values, it follows that its expectation obeys
        \begin{align*}
            \mathbb{E}[Y] = \int_0^\infty \mathrm{Pr}(\{Y \geq t\}) ~\mathrm{d}t.
        \end{align*}
        With $m' = \prod_{i=1}^n (m_i - 1)$, we obtain
        \begin{align*}
            \mathbb{E}[Y] = \int_0^\infty \mathrm{e}^{-npt} 
                \sum_{k=0}^{m'} \sum \frac{(pt)^j}{j!} \frac{j!}{i_1!i_2! \cdots i_n!} 
                ~\mathrm{d}t
                = \frac{1}{p} \sum_{j \geq 0} \sum \frac{1}{n^{j+1}} 
                    \binom{j}{i_1, i_2, \ldots, i_n}.
        \end{align*}
        The inner sum in the above equation runs over all nonnegative integers 
        $i_k$ with $k = 1,2, \ldots,n$ satisfying $i_1 + i_2 + \ldots + i_n = j$ 
        and $i_k < m_k$ for all $k$ with $1 \leq k \leq n$. The desired claim 
        follows by multiplying $\mathbb{E}[Y]$ by $p=1-q$.
    \end{pf}

    Davis and Prieditis \cite{DRPA93:spathexplength} considered the expected 
    $s$-$t$ first arrival time for the complete graph with respect to the 
    exponential model. The same analysis can be also accomplished in the case 
    of the spreading process, which is the assertion given by the following result.

    \begin{Theorem}
        \label{thm:completeGraph}
        Let $K_n = (V,E)$ be the complete graph with $n \geq 2$ vertices,
        $A \subseteq V \setminus \{t\}$, $|A|=i$ and $T_i:= T_{At}(K_n)$.
        Then $T_i$ satisfies the recurrence relation
        \begin{align*}
            T_i &= \frac{1}{1-q^{i(n-i)}} \left[
                1 + \sum_{j=i+1}^{n-1}\binom{n-1-i}{j-i} q^{i(n-j)} (1-q^i)^{j-i} T_j
            \right],
        \end{align*}
		  with the initial condition $T_{n-1}=1/(1-q^{n-1})$.
    \end{Theorem}

    \begin{pf}
        Note that in the case of the complete graph $K_n$ the transition 
        probabilities $P(A,B)$ do not depend on the sets $A$ and $B$ themselves, 
        but on their cardinality. Assume therefore that $A$ is an $i$ element 
        set and that $B$ is an $j$ element superset of $A$, then we find         
        \begin{align}
        		\label{eqn:transprobKn}
        		P(A,B) = q^{i\left( n-j\right) } \left[1-q^{i}\right]^{j-i}.
        \end{align}
				Hence, the $A$-$t$ first arrival time of $K_n$ is obtained by plugging
				\eqref{eqn:transprobKn} in the recurrence formula of the $A$-$t$ first
				arrival times stated in \eqref{eqn:momentreceqn}. That is with 
				$A \subseteq V \setminus \{t\}$ and $|A|=j$, we have
        \begin{align}
				\nonumber 
				T_{i} &= 1 + \sum_{A \subseteq C \subseteq V\setminus\{t\}} q^{i (n-|C|)}
                \left[1-q^{i}\right]^{|C|-i} ~T_{|C|}  \\
				\label{eqn:recurrenceKn}
              &= 1 + \sum_{j=i}^{n-1}\binom{n-1-i}{j-i}q^{i(n-j)}
                \left[1-q^{i}\right]^{j-i}~T_{j}.
        \end{align}
        The result then follows by factoring out $T_i$ in \eqref{eqn:recurrenceKn}.
    \end{pf}
    
    Consequently, the solution of the recurrence relation in Theorem~\ref{thm:completeGraph}
    gives us the expected $s$-$t$ first arrival time of the $K_n$ by setting $i=1$.
    In effect, the $s$-$t$ spreading resistance of $K_n$ with $n \geq 2$ yields 
    to the ratio of the $(n-1)$-st harmonic number and $n-1$. 
    This result was deduced in a different fashion in the work 
    of David and Prieditis \cite{DRPA93:spathexplength}.

    \section{Reduction techniques on series-parallel graphs}
    \label{sec:reductiontechs}
    In the case of $s$-$t$ serial parallel graphs, we present an
    algebraic approach that utilises the Hadamard product of formal 
    power series. The approach presented here is viewed as a special 
    case of Shier's \cite{SD91:netrelalgstrucs} general 
    algebraic approach for the stochastic shortest path problem
    under the assumption that all $s$-$t$ paths in $G$ are edge disjoint.

    \begin{Theorem}[{\bfseries Series reduction technique}]
        \label{thm:seriesreduction}
        Let $G = \left(V, E \right)$ be a connected graph with an articulation
        $a \in V$ such that there are two subgraphs $G_{1} = \left( V_{1}, E_{1} \right)$
        and $G_{2} = \left( V_{2}, E_{2} \right) $ of $G$ satisfying $G_{1} \cup G_{2} = G$
        and $G_{1} \cap G_{2} = \left( \left\{ a \right\} ,\emptyset \right)$. Let
        $s \in V_{1} \setminus \left\{ a \right\} $ and $t \in V_{2} \setminus
        \left\{a\right\}$, then the following equations are true:
        \begin{enumerate} [\hspace{0.25in} {\rm (}a{\rm )}]
            \item $\Phi_{st}(G,z) = \Phi_{sa}(G_1,z) \Phi_{at}(G_2,z)$
            \item $T_{st}(G) = T_{sa}(G_1) + T_{at}(G_2)$.
        \end{enumerate}
    \end{Theorem}

    \begin{pf}
        Let $\pr(\{Z_{st}(G)=n\})$ be the probability of the $s$-$t$ first arrival time with
        respect to a spread process in $G$ at $n \geq 0$ time steps. Moreover,
        denote by $\pr(\{Z_{sa}(G_1)=n\})$ and $\pr(\{Z_{at}(G_2)=n)\})$ the $s$-$a$ and $a$-$t$         
        first arrival time probabilities in $G_{1}$ and $G_{2}$ in $n\geq 0$ time steps,  
        respectively. Then we obtain
        \begin{align*}
            \pr(\{Z_{st}(G)=n\})&=\sum_{k=0}^n \pr(\{Z_{sa}(G_1)=k\})~
            \pr(\{Z_{at}(G_2)=n-k\})
        \end{align*}
        Correspondingly, let $\Phi_{st}(G,z)$ be the ordinary generating
        function for the $s$-$t$ first arrival time probability in $G$. 
        Similarly, the functions $\Phi_{sa}(G_{1},z)$ 
        and $\Phi_{at}(G_{2},z)$ denote the ordinary generating functions for the
        first arrival time that the contact process reaches $a$ in $G_{1}$ and
        $t$ in $G_{2}$, respectively. Hence
        \begin{align*}
            \Phi_{st}\left(G,z\right) &= \sum_{n \geq 0} \pr(\{Z_{st}(G)=n\})z^n
            = \Phi_{sa}\left(G_1,z\right) \Phi_{at}\left(G_2, z\right).
        \end{align*}
        follows.
        By applying the formal differentiation with respect to the
        indeterminate $z$ and then setting $z=1$ gives the desired expression
        $T_{st}(G) = T_{sa}(G_1) + T_{at}(G_2)$.
    \end{pf}

    In order to establish the parallel reduction technique  
	(Theorem~\ref{thm:parallelreduction}), we need to introduce the notion
    of Hadamard multiplication in the ring $\mathbb{C}[[Z]]$ of all formal
    power series with complex-valued coefficients over the indeterminate $z$.

    \begin{Definition}
        Let $A(z) = \sum_{i \geq 0} a_i z^i$ and $B(z) = \sum_{i \geq 0} b_i z^i$
				be arbitrary formal power series in $\mathbb{C}[[Z]]$. The \emph{Hadamard
				product} of $A(z)$ and $B(z)$ is defined as
        \begin{align*}
            (A \odot B)(z) := A(z) \odot B(z) = \sum_{i \geq 0} a_i b_i z^i.
        \end{align*}
    \end{Definition}

    \noindent In other words, the resulting formal power series $A(z) \odot B(z)$
    is obtained by termwise multiplication of $A(z)$ and $B(z)$. It is clear that the
    Hadamard product is well-defined and is closed in $\mathbb{C}[[Z]]$ as presented 
    in \cite{RS99:enumcomb}. Together with the usual polynomial addition, polynomial 
    multiplication, scalar multiplication and the Hadamard product, $\mathbb{C}[[Z]]$ 
    forms a commutative algebra over the field $\mathbb{C}$. In addition, the
    \emph{geometric series}
    \begin{align*}
        J(z) = \sum_{k \geq 0} z^k = \frac{1}{1-z}
    \end{align*}
    acts as the identity element of $\mathbb{C}[[Z]]$ under Hadamard multiplication.

		Let $r$ be a positive integer. We call the formal power series $J^r(z)$ to be 
		the \emph{$r$-geometric series}, which is defined as
		\begin{align*}
				J^r(z) = \left(\frac{1}{1-z}\right)^r = \sum_{k \geq 0} \binom{r+k-1}{k} z^k.
		\end{align*} 
		The following lemma presents a closed rational expression after pointwise
		multiplying $m$- and $n$-geometric series.
				
		\begin{Lemma}
				\label{lem:hpj}
				Let $m$ and $n$ be positive integers with $m \geq n$, then
				for any complex numbers $a$ and $b$,
				\begin{align*}
						J^m(az) \odot J^n(bz) = J^{m+n-1}(abz)
								\sum_{i=0}^{n-1} \binom{m-1}{i} \binom{n-1}{i} (abz)^i.
				\end{align*}
		\end{Lemma}
		
		\begin{pf}
				Suppose $a$ and $b$ are nonzero complex numbers. Assume 
				further that $m$ and $n$ are positive integers with $m \geq n$. 
				Let $F(z) = \sum_{k \geq 0} f(k) z^k$ be the resulting formal 
				power series obtained by taking the Hadamard product of 
				$J^m(az)$ and $J^n(bz)$. Then for each $k \in \mathbb{N}$, the 
				coefficient of $z^k$ in $F(z)$ yields
				\begin{align*}
						f(k) = \binom{m+k-1}{k} \binom{n+k-1}{k} (ab)^k.
				\end{align*}
				We see that
				\begin{align}\label{eqn:differenceFormulaEquality}
						\binom{n+k-1}{k} = \frac{(k+1)^{\overline{n-1}}}{(n-1)!}
								= \sum_{i=0}^{n-1} \binom{n-1}{i} \frac{k^{\underline{i}}}{i!}
				\end{align}
				where the rightmost part of \eqref{eqn:differenceFormulaEquality} is 
				first obtained by treating $(k+1)^{\overline{n-1}}/(n-1)!$ as a polynomial in 
				$k$ and 
				then applying the Newton's forward difference equation to the said polynomial.
				Using \eqref{eqn:differenceFormulaEquality} and setting $w = abz$ in $F(z)$, 
				we now obtain
				\begin{align}
						\label{eqn:mainthmeqn2}
						F(w) = \sum_{i=0}^{n-1} \binom{n-1}{i} \binom{m+i-1}{i} w^i ~J^{m+i}(w).
				\end{align}
				Oberve that we can factor $J^{n+m-1}(w)$ from the sum in 
				\eqref{eqn:mainthmeqn2}. This gives us
				\begin{align}
						\nonumber
						F(w) &= J^{n+m-1}(w) \sum_{i=0}^{n-1} \sum_{j=0}^{n-i-1} \binom{n-1}{i} 
								\binom{m+i-1}{i} \binom{n-i-1}{j} (-1)^j w^{i+j} \\
						\label{eqn:mainthmeqn3}
						&= J^{n+m-1}(w) \sum_{j=0}^{n-1} \sum_{i=0}^j \binom{n-1}{j-i} 
								\binom{m+j-i-1}{j-i} \binom{n-j+i-1}{i} (-1)^i w^j.
				\end{align}
				The identity
				\begin{align*}
						\binom{n-1}{j-i} \binom{m+j-i-1}{j-i} \binom{n-j+i-1}{i}
								= \binom{j}{i} \binom{m+j-i-1}{m-1} \binom{n-1}{j}
				\end{align*}
				holds for each $i$ with $0 \leq i \leq j$. In effect, $F(w)$ in
				\eqref{eqn:mainthmeqn3} becomes
				\begin{align}
						\label{eqn:mainthmeqn4}
						F(w) = J^{n+m-1}(w) \sum_{j=0}^{n-1} \binom{n-1}{j} 
								\left[ \sum_{i=0}^j \binom{j}{i} \binom{m+j-i-1}{m-1} (-1)^i \right] w^j.
				\end{align}
				
				We interpret the inner sum of \eqref{eqn:mainthmeqn4} in a combinatorial way.
				Let $Y$ be an $(m+j-1)$-element set. Let $I \subseteq Y$ be a $j$-element subset 
				of marked elements in $Y$. For each $A \subseteq I$, we define $S_A$ to be the 
				collection of $(m-1)$-element subsets of $Y$ such that all the elements of 
				$A$ are not included in the said subsets. In view of the principle of 
				inclusion-exclusion, we have
				\begin{align*}
						\left| \bigcap_{A \subseteq I} \overline{S_A} \right|
								= \sum_{A \subseteq I} (-1)^{|A|} \left| \bigcap_{a \in A} S_a \right|
								= \sum_{i=0}^j (-1)^i \binom{j}{i} \binom{m+j-i-1}{m-1}. 
				\end{align*} 
				This gives us the number of ways we can draw an $(m-1)$-element set 
				from $Y$ such that all the marked elements of $I$ are included. In turn, 
				this corresponds to the number of ways of forming an $j$-element subset 
				of an arbitrary $(m-1)$-element set.
		\end{pf}

    \begin{Theorem}[{\bfseries Parallel reduction technique}]
        \label{thm:parallelreduction}
        Let $G = \left( V, E \right)$ be a graph and $\left\{ s, t \right \}$
        be a separating vertex pair of $G$ such that $H$ and $K$ are two connected
        nontrivial graphs with $H \cup K = G$ and $H \cap K= \left( \left\{ s, t \right\},
        \emptyset \right)$, then
        \begin{align*}
            \Phi_{st}(G,z) &= \Phi_{st}(H,z) + \Phi_{st}(K,z)
                -(1-z) \left[ \frac{\Phi_{st}(H,z)}{1-z} \odot \frac{\Phi_{st}(K,z)}{1-z} \right]
        \end{align*}
    \end{Theorem}

    \begin{pf}
        Let $Z_{st}(G)$ be the random variable of the $s$-$t$ first arrival time in $G$, 
        and $Z_{st}(H)$ and $Z_{st}(K)$ be the random variables of the $s$-$t$ first arrival
        times in $H$ and $K$, respectively. As the
        events $\{Z_{st}(H) > n\}$ and $\{Z_{st}(K) > n\}$ are independent, we have
        \begin{align*}
            1 - \pr\left( \{Z_{st}(G) \le n\} \right) = 
            		\left( 1 - \pr\left( \{Z_{st}(H) \le n\} \right) \right)
                \left( 1 - \pr \left( \{ Z_{st}(K) \le n \} \right) \right).
        \end{align*}
        Multiplying both sides of the above equation with $z^n$ and
        summation over all $n \ge 0$ yields
        \begin{align*}
            \frac{1 - \Phi_{st}(G,z)}{1 - z} &=
            \frac{1}{1 - z} - \frac{\Phi_{st}(H,z)}{1 - z} - \frac{\Phi_{st}(K,z)}{1 - z}
            + \frac{\Phi_{st}(H,z)}{1 - z} \odot \frac{\Phi_{st}(K,z)}{1 - z},
        \end{align*}
        which leads after rearranging terms to the desired result.
    \end{pf}

    We apply the parallel reduction for the series-parallel graph $G$
    with $s$ and $t$ as its start and terminal vertices such that $\{s,t\}$
    forms a vertex separating pair. The subsequent result deals with the 
    case that one of the two subgraphs of $G$ is a path of length $r \geq 1$.

    \begin{Corollary}
        Let $G = \left( V,E \right)$ be a graph and $\{s,t\}$ a separating
        vertex pair of $G$, furthermore assume that $H$ and $K$ are subgraphs 
        of $G$ with $G = H \cup K$ and $H \cap K = (\{s,t\}, \emptyset)$ and 
        that $K$ is a path of length $r$ connecting $s$ and $t$. Then the 
        following equations hold:
        \begin{enumerate}[\hspace{0.25in} {\rm (}a{\rm )}]
            \item $\Phi_{st}(G,z) = \displaystyle \frac{p^{r}z^{r}}{(1-qz)^{r}} +
                \left( \frac{p}{q} \right)^{r-1} \frac{1-z}{(r-1)!}
                \frac{d^{r-1}}{dz^{r-1}}\left( \frac{z^{r-1}\Phi_{st}(H,qz)}{1-qz} \right)$
            \item $T_{st}(G) = \displaystyle \frac{rq^r}{(1-q)^{r+1}}-\frac{1}{(r-1)!}
                \left(\frac{p}{q}\right)^{r-1} \frac{d^{r-1}}{dz^{r-1}} \left.
                    \left(\frac{z^{r-1}\Phi_{st}(H,qz)}{1-qz}\right)
                \right\vert_{z=1}$.
        \end{enumerate}
    \end{Corollary}

    \begin{pf}
        For the first arrival times of the path $K$ of length $r$ in $G$, one finds
        \begin{align}
            \label{eqn:ogffatpath}
            \Phi_{st}(K,z) = \sum_{n\geq r} \binom{n-1}{r-1}p^{r}q^{n-r}z^{n}
                = \frac{p^{r}z^{r}}{(1-qz)^{r}}.
        \end{align}
        By Theorem \ref{thm:parallelreduction}, the ordinary generating function for the
        first arrival time of $G$ yields
        \begin{align*}
            \Phi_{st}\left(G,z\right) &= \Phi_{st}\left(H,z\right) + \frac{p^{r}z^{r}}{(1-qz)^{r}} +
                \left(z-1\right) \left[ \frac{\Phi_{st}\left(H,z\right)}{1-z} \odot
                \frac{p^{r}z^{r}}{\left( 1-z\right) (1-qz)^{r}}\right] \\
            &= \Phi_{st}\left(H,z\right) + \frac{p^{r}z^{r}}{(1-qz)^{r}} + \left( z-1\right)
                \left[\frac{\Phi_{st}\left(H,z\right) }{1-z} \odot
                \left( \frac{1}{1-z} - \left( \frac{p}{q} \right)^{r-1}
                    \frac{1}{\left( 1 - qz \right)^{r}}\right) \right] \\
            &= \frac{p^{r}z^{r}}{(1-qz)^{r}} + \left( \frac{p}{q}\right)^{r-1}
                \left(1-z\right) \left[ \frac{\Phi_{st}\left(H,z\right)}{1-z} \odot
                \frac{1}{\left( 1-qz\right) ^{r}}\right] \\
            &= \frac{p^{r}z^{r}}{(1-qz)^{r}} + \left( \frac{p}{q}\right)^{r-1}
                \left( 1 - z \right) \left[ \frac{1}{\left( 1-z\right) ^{r}} \odot
                \frac{\Phi_{st} \left(H,qz\right)}{1-qz} \right] \\
            &= \frac{p^{r}z^{r}}{(1-qz)^{r}} + \left( \frac{p}{q}\right)^{r-1}\frac{1-z}{(r-1)!}
                \frac{d^{r-1}}{dz^{r-1}}\left(\frac{z^{r-1}\Phi_{st}\left(H,qz\right)}{1-qz}\right).
        \end{align*}
        After evaluating the formal derivative at $z = 1$ we obtain
        \begin{align*}
            &T_{st}(G) = \frac{rq^r}{(1-q)^{r+1}} - \frac{1}{(r-1)!}
                \left(\frac{p}{q}\right)^{r-1} \frac{d^{r-1}}{dz^{r-1}}
                \left.
                    \left( \frac{z^{r-1}\Phi_{st}(H,qz)}{1-qz} \right)
                \right\vert_{z=1}. \qedhere
        \end{align*}
    \end{pf}

    Another instance of the parallel reduction technique is given by two parallel paths
    of lengths $n$ and $m$ between $s$ and $t$.

    \begin{Corollary}
        \label{cor:twopaths}
        Let $G=(V,E)$ be the graph that consists of two parallel paths of lengths 
        $n \geq 1$ and $m \geq 1$ between the vertices $s,t\in V$ where $m \geq n$. Then
        \begin{align*}
            \Phi_{st}(G,z) &= 1 - \left(\frac{p}{q}\right)^{n+m-2} \frac{1-z}{(1-q^2z)^{m+n-1}}
                \sum_{\ell=0}^{n-1} \binom{n-1}{\ell} \binom{m-1}{\ell} (q^2z)^\ell.
        \end{align*}
    \end{Corollary}

    \begin{pf}
        Let $\Phi_{st}(P_\ell,z)$ be the ordinary generating function
        for the $s$-$t$ first arrival time in a path of length $\ell$.
        Then the ratio $\Phi_{st}(P_\ell,z) / (1-z)$ becomes
        \begin{align*}
            \frac{\Phi_{st}(P_\ell,z)}{1-z} &= \frac{1}{1-z} -
                \left( \frac{p}{q} \right)^{\ell-1} \frac{1}{(1-qz)^\ell}
        \end{align*}
        by Equation~\eqref{eqn:ogffatpath} and the application of partial fraction
        decomposition. With $\ell = n$ and $\ell = m$, it follows from the
        parallel reduction technique that
        \begin{align}
            \label{eqn:hprodjfuncs}
            \Phi_{st}(G,z) &= 1 - (1-z) \left( \frac{p}{q} \right)^{n+m-2}
                \left[ \frac{1}{(1-z)^n} \odot \frac{1}{(1-q^2z)^m} \right].
        \end{align}
        The desired result now follows by applying Lemma \ref{lem:hpj}
        to \eqref{eqn:hprodjfuncs}.
    \end{pf}

    \section{Bounds for the expected first arrival time}
		\label{sec:bounds}
    As a result of the exponentially increasing number of equations in
    Corollary \ref{cor:momentreceqn}, the exact computation of the expected
    $s$-$t$ first arrival times can be performed in reasonable time only in small
    networks. Hence, computational feasible bounds for $T_{st}(G)$ are of great interest.

    \begin{Theorem}
        Let $G=(V,E)$ be a graph with $s,t\in V$ and $H=(V,F)$ with $F \subseteq E$
        being a subgraph of $G$, then $T_{st}(H) \geq T_{st}(G)$.
    \end{Theorem}

    \begin{pf}
        Every removal of an edge $e\in E$ of $G$ causes an increase in the
        expected $s$-$t$ first arrival time.
    \end{pf}

    By considering a shortest $s$-$t$ path as a subgraph of $G$ one gets the
    following corollary giving an upper bound, that is not very tight in general.

    \begin{Corollary}
        Let $G=(V,E)$ be a graph and $s,t\in V$. If $d(s,t)$ is the 
        distance of $s$ and $t$ in $G$, then the inequality
        \begin{align*}
        	\mathbb{E}[T_{st}(G)] \leq \dfrac{d(s,t)}{1-q}
        \end{align*}
        always holds. \hfill $\blacksquare$
    \end{Corollary}

    \noindent A known lower bound for $\tau_{st}(G)$ was presented by Lyons
    \cite{LRPRPY99:resistancebound} and is connected to the ideas presented in 
    \cite{NWC59:randwalkeleccurrents}. That is, if $\mathrm{Res}_{st}(G)$
    is the electrical resistance between the vertices $s$ and $t$ in a
    given connected graph $G$, where the edges are assumed to have unit resistance, then
    \begin{align*}
        \frac{\mathrm{Res}_{st}(G)}{1-q} \leq \tau_{st}(G).
    \end{align*}





    \noindent A lower bound for $T_{st}(G)$ can be found by considering 
    the $s$-$t$ reliability polynomial of $G$.
    
    \begin{Definition}
        Let $G = (V,E)$ be a connected graph with $m$ edges and
        $R_{st}(G,q) = c_0 + c_1q + c_2q^2 + \ldots + c_mq^m$
        the $s$-$t$ reliability polynomial of $G$, i.e. the probability that there is at
        least one intact $s$-$t$ path in $G$, if the edges of $G$ are failing independently with  
        probability q. The 
        \emph{$s$-$t$ insertion probability}, denoted by 
        $\widetilde{R}_{st}(G,q)$, is then defined as
        \begin{align*}
            \widetilde{R}_{st}(G,q) &= \sum_{i=1}^m \frac{c_i}{q^i-1}.
        \end{align*}
    \end{Definition}
    
    \begin{Theorem}
        \label{thm:sophisticatedRelbound}
        Let $G=(V,E)$ be a graph and $s,t\in V$. If $d(s,t)$ is the distance between $s$ and 
		  $t$ in $G$, then
        \begin{align*}
            d(s,t) - 1 + \widetilde{R}_{st}(G,q)\le T_{st}(G).
        \end{align*}
    \end{Theorem}

    \begin{pf}
        For each edge $e\in E$, let $L(e)$ be independent geometric random variables, with
        non infection probabilities $q$, representing the edge lengths. Denote
        now by $G_k=(V, E_k)$ with $k=0,1,2,\ldots$ a discrete random process that 
        gives an infinite sequence of subgraphs $G_k$ of $G$, where $e\in E_k$ 
        if and only if the event $\{L(e) \leq k\}$ occurs. 
        Moreover, denote by $I_{st}$ the random variable
        of the smallest $k$, so that $G_k$ contains an $s$-$t$ path.

        Now the following inequality holds
        \begin{align}
            \label{eqn:fundamentalInequality}
            \mathrm{Pr}\left(\{I_{st} \leq k\}\right) \geq
                \mathrm{Pr}\left(\{Z_{st} \leq d(s,t)+k-1\}\right),
        \end{align}
        for all $k \geq 0$. In order to prove this estimate, we assume an 
        elementary event $E$ in \linebreak $\{Z_{st} \leq d(s,t)+k-1\}$, e.g. a
        realisation of edge lengths, so that the length of a shortest
        $s$-$t$ path is smaller or equal than $d(s,t)+k-1$. In order
        that $E$ occurs there must be at least one $s$-$t$ path
        in $G$, so that the sum of the edge lengths of the path is
        smaller or equal than $d(s,t)+k-1$. Therefore all edge lengths 
        of this $s$-$t$ path must be smaller than $k$, as if there would be one
        edge length in this path that is greater than $k+1$, then the sum of the path
        lengths would exceed $(k+1)+(d(s,t)-1)$, which contradicts the assumption
        $Z_{st} \leq d(s,t)+k-1$. As all edge lengths are smaller than $k$, the event $E$ 
        is by definition also in $\{I_{st}\le k\}$, i.e. $G_k$ must also contains this 
        $s$-$t$ path. Furthermore the probabilities for the event $E$ coincide in both  
        probability spaces, so that \eqref{eqn:fundamentalInequality} follows.
        
        Due to \eqref{eqn:fundamentalInequality}, we find
        that the expectation of the random variable $I_{st}$ is smaller
        than the expectation of the random variable $Z_{st}-d(s,t)+1$. Hence,
        \begin{align*}
             \mathbb{E}\left[I_{st}\right] \leq \mathbb{E}\left[Z_{st}-d(s,t)+1\right]= 
            \mathbb{E}\left[Z_{st}\right]-d(s,t)+1
        \end{align*}
        and this gives $d(s,t)-1+\mathbb{E}\left[I_{st}\right]$ as a lower bound for $T_{st}(G)$.

        In order to show the relationship between $\mathbb{E}[I_{st}]$ and the 
        $s$-$t$ reliability polynomial,
        consider the edge $e = \{u,v\}$ and the event $\{L(e)\le k\}$.
        We easily find that $\mathrm{Pr}(\{L(e)\le k\}) = 1-q^k$, whenever
        the non insertion probability is $q$. Assume now that there are
        $k$ parallel edges between $u$ and $v$, that are failing independently
        with probability $q$, then the probability that there is at least one non
        failing edge is also given by $1-q^k$. Therefore, the probability for the
        event $\{I_{st} > k\}$ coincides with the probability of the event that
        there is no $s$-$t$ path in the multigraph $G(k)$, where $G(k)$ emerges from
        $G$ by replacing each edge in $G$ with $k$ parallel edges. Equivalently,
        \begin{align*}
            \pr(\{I_{st} > k\}) = 1-R_{st}(G(k), q) = 1-R_{st}(G,q^k)
        \end{align*}
        and we can write
        \begin{align*}
            \mathbb{E}\left[I_{st}\right] = \sum_{k \ge 1} k \left[R_{st}(G,q^k,\ldots, q^k)
                - R_{st}(G, q^{k-1},\ldots, q^{k-1}) \right].
        \end{align*}
        Now suppose that the $s$-$t$ reliability polynomial is
        expressed as $R_{st}(G, q) = \sum_{i=0}^m c_i q^i$, then
        \begin{align*}
            \mathbb{E}\left[I_{st}\right]
                = \sum_{i=0}^m c_i \sum_{k\ge 1} k ~(q^{ik}-q^{i(k-1)})
            = \sum_{i=1}^m \frac{c_i}{q^i-1} = \widetilde{R}_{st}(G,q)
        \end{align*}
        and hence, $d(s,t)-1 + \widetilde{R}_{st}(G,q) \leq T_{st}(G)$.
    \end{pf}

    \section{Conclusion}
    We have shown that there is a one-to-one correspondence of the
    spread process and the stochastic shortest path problem.
    Furthermore, it was shown that this process can be well-described
    by the exponential model, when the infection probability $p$ tends
    to zero. Several closed formulae or at least efficient calculation 
    schemes for special graphs, such as complete graphs, parallel paths
    and $s$-$t$ series-parallel graphs are presented. Finally some ideas 
    for bounds are proposed that utilise the well-known concept of the 
    $s$-$t$ reliability polynomial, which can be efficiently computed for graphs 
    with bounded tree-width.

    The authors conjecture the following bound for the $s$-$t$ 
    first arrival time, that behaves good for values of $q$, that are close to $1$, i.e. 
    $\lim_{q\rightarrow 1} T_{st}(G)/\tau_{st}(G)=1$.
    
    \begin{Conjecture}
        \label{conj:lowerBoundLimitTheorem}
        Let $G=(V,E)$ be a graph, $s,t \in V$ and $T_{st}(G)$,
        $\tau_{st}(G)$ be the expected $s$-$t$ first arrival times of the
        spread process and the exponential model, respectively. In this case
        $\tau_{st}(G) \leq T_{st}(G)$ holds for all $0 \leq q < 1$.
    \end{Conjecture}


    \section{Acknowledgements}
    The authors Raymond Lapus and Frank Simon receive a grant from the 
    German Academic Exchange Service (DAAD) and the European Social Fund (ESF), 
    respectively. The authors want to thank Frank G\"{o}ring from the 
    University of Technology in Chemnitz for giving an direct and nice symmetry proof 
    for Corollary~\ref{cor:symmetry}.

    \bibliographystyle{plain}
    \bibliography{references}
\end{document}